\documentclass[12pt,twoside, A4paper]{article}
\usepackage[papersize={21cm,29.7cm},left=3.5cm,top=3cm,right=3.5cm,bottom=2.5cm]{geometry}
\usepackage{latexsym,enumerate}
\usepackage{amsmath,amsthm,amsopn,amstext,amscd,amsfonts,amssymb}
\usepackage{graphics,graphicx}
\usepackage{setspace}
\usepackage[english]{babel}
\usepackage[T1]{fontenc}

\usepackage{blindtext}
\usepackage{titlesec}
\singlespace

\newtheorem{theorem}{Theorem}[section]
\newtheorem{definition}{Definition}[section]

\newtheorem{lem}{Lemma}[section]

\newtheorem{rem}{Remark}[section]

\def\title#1{\vspace{24pt}{ \bf \begin{center} \hrule\vspace{60pt}\fontsize{12pt}{12pt}\selectfont Lecture / Ponencia XX \\ \vspace{14pt} \fontsize{14pt}{14pt}\selectfont #1 \vspace{0pt} \end{center}  } }
\def\authors#1{{ \begin{center} \fontsize{12pt}{12pt}\selectfont  #1 \vspace{0pt} \end{center} } }
\def\affiliation#1{{\sl \begin{center}\vspace{10pt} \fontsize{10pt}{10pt}\selectfont #1 \vspace{0pt} \end{center} } }
\def\name#1{\unskip$^{#1}$}

\titleformat{\section}[hang]
{\bf}{\thesection.}{0.7cm} {}

\fontsize{12pt}{12pt}

\begin{document}


\title{Hitting times of threshold exceedances and their distributions}

\authors{N. Markovich\name{1}}

\affiliation{\name{1} Institute of Control Sciences,
                Russian Academy of Sciences, Russia,
nat.markovich@gmail.com, markovic@ipu.rssi.ru}




\section*{Abstract}


We investigate exceedances of the process over a sufficiently high threshold. The exceedances determine the risk of hazardous events like climate catastrophes, huge insurance claims, the loss and delay in telecommunication networks.
 Due to dependence such exceedances tend to occur in clusters. Cluster structure of social networks is caused by dependence (social relationships and
interests) between nodes and possibly heavy-tailed distributions of the node
degrees. A minimal time to reach a large node determines the first hitting time.
We derive asymptotically equivalent distribution and a limit expectation of the first
hitting time to exceed the threshold $u_n$ as sample size $n$ tends to infinity. The results can be extended to the second and, generally, to $k$th ($k>2$) hitting times.

\section{Introduction}

Let $\{X_n\}_{n\ge 1}$ be a stationary process with marginal distribution function $F(x)$ and $M_n=\max\{X_1,...,X_n\}$. We investigate the exceedances of the process over a sufficiently high threshold $u$. Due to dependence such exceedances tend to occur in clusters.
 Let us consider the inter-cluster size
 \begin{eqnarray}\label{18}T_1(u)=\min\{j\ge 1: M_{1,j}\le u,X_{j+1}>u|X_{1}>u\},\end{eqnarray} i.e. the number of inter-arrivals of observations running under the threshold between two consecutive exceedances,
where $M_{1,j}=\max\{X_{2},...,X_{j}\}$, $M_{1,1}=-\infty$.
\\
Let $T^*(u_n)$ be a first hitting time to exceed the threshold $u_n$. We get
\begin{eqnarray}\label{0}P\{T^*(u_n)=j+1\}&=& P\{M_j\le u_n, X_{j+1}>u_n\},\end{eqnarray} $j=0,1,2,...$, $M_{0}=-\infty$.
\\
The necessity to evaluate quantiles of the  first hitting time and its mean is arisen in many applications. In social networks it is important to compare sampling strategies (Avrachenkov \textsl{et al.} (2012), Avrachenkov \textsl{et al.} (2014a), Avrachenkov \textsl{et al.} (2014b), Chul-Ho Lee \textsl{et al.} (2012)) like  random walks, Metropolis-Hastings Markov chains, Page Ranks and other with regard to how quickly they allow to reach a node with a large degree that is the number of links with other nodes. It is important to investigate the first hitting time of significant nodes since it allows us to disseminate advertisement or collect opinions more effectively within the clusters surrounded such nodes. Similar problem is required in telecommunication peer-to-peer networks to find a node with a large number of peers, D\'{a}n and  Fodor (2009), Markovich (2013).
\begin{definition} (Leadbetter \textsl{et al.} (1983), p.53) The stationary sequence  $\{X_n\}_{n\ge 1}$ is said to have extremal index
$\theta\in[0,1]$ if
for each $0<\tau <\infty$ there is a sequence of real numbers $u_n=u_n(\tau)$ such that
\begin{equation}\label{1}\lim_{n\to\infty}n(1-F(u_n))=\tau \qquad\mbox{and}\end{equation}
\begin{equation}\label{2}\lim_{n\to\infty}P\{M_n\le u_n\}=e^{-\tau\theta}\end{equation}
hold.
\end{definition}
The  extremal index $\theta$ of  $\{X_n\}$ relates to the first hitting time  $T^*(u_n)$, Roberts \textsl{et al.} (2006).
Really, since
 $u_n$ is selected according to (\ref{1}) it follows that $P\{X_n > u_n\}$ is asymptotically equivalent to $1/n$. Notice, that $P\{M_k\le u_n\}=P\{T^*(u_n)>k\}$. Hence, substituting $\tau$ by (\ref{1})  we get from (\ref{2}) \[P\{T^*(u_n)/n>k/n\}\sim e^{-\theta k P\{X_n > u_n\}}\sim e^{-\theta k/n}\] and
 \[\lim_{n\to\infty}P(T^*(u_n)/n > x) = e^{-\theta x}\] for positive $x$. It follows
 \begin{equation}\label{5b}\lim_{n\to\infty}E(T^*(u_n)/n)=1/\theta.\end{equation}
 This implies,  the smaller $\theta$, the longer it takes to reach an observation with a large value.
\\
Using achievements regarding the asymptotically equivalent geometric distribution of $T_1(x_{\rho_n})$ derived in (Theorem \ref{T2}, Markovich (2014)), where the $(1-\rho_n)$th quantile $x_{\rho_n}$ of $\{X_n\}$ is taken as $u_n$, we derive an asymptotically equivalent distribution of the first hitting time and its limiting expectation that specifies (\ref{5b}) in Section \ref{Sec2}.
\\
We use the mixing condition proposed in Ferro and Segers (2003)
\begin{equation}\label{44}\alpha_{n,q}(u)=\max_{1\le k\le n-q}\sup|P(B|A)-P(B)|=o(1),\qquad n\to\infty,
\end{equation}where for real  $u$ and integers $1\le k\le l$,
$\mathcal{F}_{k,l}(u)$
is the $\sigma$-field generated by the events $\{X_i>u\}$, $k\le
i\le l$ and the supremum is taken over all $A\in
\mathcal{F}_{1,k}(u)$ with $P(A)> 0$ and $B\in
\mathcal{F}_{k+q,n}(u)$ and $k$, $q$ are positive integers.
\\
We use the following partition \begin{equation}\label{31}
1\le k_{n,1}^*\le k_{n,2}^*\le k_{n,3}^*\le  k_{n,4}^*\le j,\qquad
j=2,3,...,
\end{equation} $k_{n,0}^*=1$, $k_{n,5}^*=j$,
$k_{n,i}^*=[jk_{n,i}/n]+1, i=\{1,2\}$, $k_{n,3}^*=j-[jk_{n,4}/n]$, $k_{n,4}^*=j-[jk_{n,3}/n]$,
\footnote{$[x]$ represents the integer part of the real number $x$.}  of the interval $[1,j]$ for a fixed $j$, where positive integers $\{k_{n,i}\}$ are such that \begin{equation}\label{4a}\{k_{n,i-1}=o(k_{n,i}), i\in\{2,3,4\}\}, \qquad  k_{n,4}=o(n).\end{equation}
\begin{theorem}\label{T2} (Markovich (2014))
Let  $\{X_n\}_{n\ge 1}$ be a stationary process with the extremal index 
$\theta$.
Let $\{x_{\rho_n}\}$ and $\{x_{\rho_n^*}\}$ be  sequences of quantiles of $X_1$ of the levels $\{1-\rho_n\}$ and $\{1-\rho_n^*\}$, respectively,\footnote{$\overline{F}(x_{\rho_n})=P\{X_1>x_{\rho_n}\}=\rho_n$.
}
those satisfy the conditions (\ref{1}) and (\ref{2}) if $u_n$ is replaced by $x_{\rho_n}$ or by $x_{\rho_n^*}$  and, $q_n=1-\rho_n$, $q_n^*=1-\rho_n^*$, $\rho_n^*=(1-q_n^{\theta})^{1/\theta}$.
Let positive integers $\{k^*_{n,i}\}$, $i=\overline{0,5}$, and $\{k_{n,i}\}$, $i=\overline{1,4}$, be  respectively as in  (\ref{31}) and (\ref{4a}),
$p_{n,i}^*=o(\Delta_{n,i})$, $\Delta_{n,i}=k_{n,i}^*-k_{n,i-1}^*$, $q_{n,i}^*=o(p_{n,i}^*)$, $i\in\{1,2,...,5\}$ and $\{p_{n,3}^*\}$ be an increasing sequence,
such that \begin{eqnarray}\label{29}
\alpha^*_n(x_{\rho_n})&=&
\max\{\alpha_{k_{n,4}^*,q_{n,1}^*}; \alpha_{k_{n,3}^*,q_{n,2}^*}; \alpha_{\Delta_{n,3},q_{n,3}^*}; \alpha_{j+1-k_{n,2}^*,q_{n,4}^*};
\\
&&
\alpha_{j+1-k_{n,1}^*,q_{n,5}^*};  \alpha_{j+1,k_{n,4}^*-k_{n,1}^*}
\}=o(1)\nonumber\end{eqnarray}
and
\begin{eqnarray}\label{29a} \alpha_{j+1,k_{n,4}^*-k_{n,1}^*}&=& o(\rho_n)
\end{eqnarray}
hold as $n\to\infty$,
 where $\alpha_{n,q}=\alpha_{n,q}(x_{\rho_n})$ is determined by (\ref{44}).
 Then  it holds for $j\ge 2$
\begin{eqnarray}\label{2a}
\lim_{n\to\infty}P\{T_1(x_{\rho_n})=j\}/(\rho_n(1-\rho_n)^{(j-1)\theta})&=&1,
\end{eqnarray}
\end{theorem}
The achievements can be similarly extended to the second and, generally, to $k$th ($k>2$) hitting times.
The paper is organized as follows. In Section \ref{Sec2} we derive the asymptotically equivalent distribution and the limit expectation of the  first hitting time to exceed sufficiently high threshold. Proofs are given in Section \ref{Proof}.

\section{Distribution and expectation of the first hitting time}\label{Sec2}
\begin{theorem}\label{Theor1} Let all conditions of Theorem \ref{T2} are satisfied. We assume
\begin{equation}\label{3e}
\sup_n \rho_nE (T_1(x_{\rho_n}))<\infty.
\end{equation}
Then we get
\begin{eqnarray}\label{19}\lim_{n\to\infty}nP\{T^*(u_n)=n\}&=& e^{-\theta\tau}/\theta,\end{eqnarray}
 \begin{eqnarray}\label{17}
P\{T^*(x_{\rho_n})=j\}&\sim & \frac{\rho_n}{\theta}(1-\theta\rho_n)^{j-1}.
\end{eqnarray}
as $n\to\infty$.
 \end{theorem}
Expression (\ref{17}) implies that $P\{T^*(x_{\rho_n})=j\}\sim  P\{T_{\theta}=j\}$, where
\begin{eqnarray*}
T_{\theta}&=&\left\{\begin{array}{ll}
                                       \chi, & \mbox{with probability}\qquad 1/\theta^2; \\
                                        0, & \mbox{with probability}\qquad 1-1/\theta^2
                                      \end{array}
                                    \right.
\end{eqnarray*}
holds and $\chi$ has the geometric distribution with probability $\rho_n\theta$.
\begin{rem}\label{Rem6}
The condition (\ref{3e})
provides a uniform convergence of the range $\sum_{j=1}^{\infty}P\{T_1(x_{\rho_n})=j\}$ by $n$.
 The latter condition is fulfilled for a geometrically distributed $T_1(x_{\rho_n})$, i.e. $P\{T_1(x_{\rho_n})=j\}=\rho_n(1-\rho_n)^{j-1}$.
 \end{rem}
\begin{lem}\label{Lem1}Let conditions of Theorem \ref{Theor1} be satisfied.
Then it follows
\begin{equation}\label{5a}
\lim_{n\to\infty}\rho_nET^*(x_{\rho_n})=1/\theta^3.
\end{equation}
\end{lem}
Since $\rho_n\sim \tau/n$ according to (\ref{1}), the expression (\ref{5a}) specifies (\ref{5b}). 
\section{Proofs}\label{Proof}
\subsection{Proof of Theorem \ref{Theor1}}
It follows from (\ref{0}) that
\begin{eqnarray}\label{3}P\{T^*(u_n)=n\} &=& P\{M_{n-1}\le u_n, X_{n}>u_n\}\nonumber
\\
&=& P\{M_{n-1}\le u_n\}-P\{M_{n}\le u_n\}.
\end{eqnarray}
We obtain due to (\ref{18})  and the stationarity of $\{X_n\}$  that it holds
\begin{eqnarray}\label{13}&& P\{T_1(u_n)=n\} = P\{M_{1,n}\le u_n, X_{n+1}>u_n|X_{1}>u_n\}\nonumber
\\
&=& \left(P\{M_{1,n}\le u_n, X_{n+1}>u_n\}-P\{M_{n}\le u_n, X_{n+1}>u_n\}\right)/P\{X_{1}>u_n\}\nonumber
\\
&=& \left(P\{M_{n-1}\le u_n, X_{n+1}>u_n\}-P\{M_{n}\le u_n, X_{n+1}>u_n\}\right)/P\{X_{1}>u_n\}\nonumber
\\
&=&\left(P\{T^*(u_n)=n\}-P\{T^*(u_n)=n+1\}\right)/P\{X_{1}>u_n\}.
\end{eqnarray}
Following Ferro and Segers (2003)  we get alternatively for $n\ge 1$
\begin{eqnarray*}P\{T_1(u_n)>n\} &=& P\{M_{1,n+1}\le u_n| X_{1}>u_n\}
\\
&=&\left(P\{M_{1,n+1}\le u_n\}-P\{M_{n+1}\le u_n\}\right)/P\{X_{1}>u_n\}
\\
&=&\left(P\{M_{n}\le u_n\}-P\{M_{n+1}\le u_n\}\right)/P\{X_{1}>u_n\}
\\
&=& P\{T^*(u_n)=n+1\}/P\{X_{1}>u_n\}.
\end{eqnarray*}
Thus, we get
\begin{eqnarray}\label{4}
P\{T^*(x_{\rho_n})=n+1\}&=& P\{X_{1}>x_{\rho_n}\}\cdot P\{T_1(x_{\rho_n})>n\}\nonumber
\\
&=& \rho_n\sum_{i=n+1}^{\infty}P\{T_1(x_{\rho_n})=i\},
\end{eqnarray}
where $\rho_n= P\{X_{1}>x_{\rho_n}\}$. Now we use the fact followed from (\ref{2a}) in Theorem \ref{T2}, i.e. it holds
\begin{eqnarray}\label{14}c_nP\{T_1(x_{\rho_n})=i\}\sim \eta_n(1-\eta_n)^{i-1}\end{eqnarray} as $n\to\infty$, where $1-\eta_n=(1-\rho_n)^{\theta}$, $c_n=\eta_n/(1-(1-\eta_n)^{1/\theta})$, $0<\eta_n<1$.
\\
Furthermore, we obtain
\begin{eqnarray}\label{5}
P\{T^*(x_{\rho_n})=n+1\}&=& \frac{\rho_n}{c_n}\sum_{i=n+1}^{\infty}c_nP\{T_1(u_n)=i\} \sim \frac{\rho_n}{c_n}\sum_{i=n+1}^{\infty}\eta_n(1-\eta_n)^{i-1}\nonumber
\\
&=&\frac{\rho_n}{c_n}(1-\eta_n)^{n}= \frac{\rho_n}{\eta_n}\left(1-(1-\eta_n)^{1/\theta}\right)(1-\eta_n)^{n}\nonumber
\\
&\sim &\frac{\rho_n^2(1-\theta\rho_n)^{n}}{1-(1-\rho_n\theta)}=\frac{\rho_n}{\theta}(1-\theta\rho_n)^{n}
\sim \frac{\rho_n}{\theta}e^{-\theta\tau},
\end{eqnarray}
since $(1-\rho_n)^{\theta}=1-\theta\rho_n+o(\rho_n)$ and $\lim_{n\to\infty}(1-\theta\rho_n)^{n}=e^{-\theta\tau}$ as $\rho_n\sim\tau/n$, $n\to\infty$. Formula (\ref{5}) implies both (\ref{19}) rewriting it with regard to $P\{T^*(x_{\rho_n})=n+1\}/\rho_n$  and (\ref{17}) replacing $n+1$ by $j$.
\\
Hence, we have from (\ref{5})
\begin{eqnarray}\label{6}
P\{T^*(x_{\rho_n})=j\}&\sim & \frac{1}{\theta^2}\rho_n\theta(1-\rho_n\theta)^{j-1}, 
\qquad n\to\infty.
\end{eqnarray}
\\
The similarity in the first string of (\ref{5}) follows from (\ref{3e}) and  since
\begin{eqnarray*}\sup_n r_k(n)&=& \sup_n \sum_{i=\lfloor k/\rho_n\rfloor}^{\infty}P\{T_1(x_{\rho_n})=i\}
\\
&=&
\sup_n \sum_{i=\lfloor k/\rho_n \rfloor}^{\infty}
\frac{i}{i}P\{T_1(x_{\rho_n})=i\}\le \sup_n \frac{\rho_n}{k}E (T_1(x_{\rho_n}))\to 0
\end{eqnarray*}
holds as $k\to\infty$.

\subsection{Proof of Lemma \ref{Lem1}}
Using (\ref{4}) and (\ref{14}) we find the expectation of the first hitting time
\begin{eqnarray*}
ET^*(x_{\rho_n})&=&\sum_{j=1}^{\infty}jP\{T^*(x_{\rho_n})=j\}=\sum_{j=1}^{\infty}j\rho_n\sum_{i=j}^{\infty}P\{T_1(x_{\rho_n})=i\}
\\
&\sim &  \sum_{j=1}^{\infty}j\frac{\rho_n}{c_n}\sum_{i=j}^{\infty}\eta_n(1-\eta_n)^{i-1}=\frac{\rho_n}{c_n\eta_n}\sum_{j=1}^{\infty}j\eta_n(1-\eta_n)^{j-1}
\\
&=& \frac{\rho_n^2}{\eta_n^3}=\frac{\rho_n^2}{(1-(1-\rho_n)^{\theta})^3}.
\end{eqnarray*}
The similarity follows from the same arguments as before. Then (\ref{5a}) follows.




\section*{Acknowledgments}

Special thanks to the partial financial support by  the Russian Foundation for Basic Research, grant 13-08-00744 A.



\section*{References}

\begin{description}
\item Avrachenkov, K., Litvak, N., Sokol, M. and Towsley, D. (2012).
{\it Algorithms and Models for the Web Graph}.   Lecture Notes in Computer Science, {\bf 7323},
Quick Detection of Nodes with Large Degrees, Springer Berlin Heidelberg, 54-65.
\item Avrachenkov, K.,  Markovich, N. and Sreedharan, J. K. (2014a). Distribution and Dependence of Extremes in Network Sampling Processes.
   {\it Research report no.\ 8578, INRIA, http://hal.inria.fr/hal-01054929}
\item Avrachenkov, K.,  Markovich, N. and Sreedharan, J. K. (2014b). Distribution and Dependence of Extremes in Network Sampling Processes. Workshop on Complex Networks and their Applications. {\it The 10th IEEE International Conference on Signal-Image Technology and Internet-Based Systems (SITIS 2014) November 23-27, 2014, Marrakech, Morocco}, 331-338.
\item Chul-Ho Lee, Xin Xu and Do Young Eun. (2012).	Beyond Random Walk and Metropolis-Hastings Samplers: Why You Should Not Backtrack for Unbiased Graph Sampling. 	{\it CoRR}.
\item  D\'{a}n, G. and  Fodor, V. (2009). Delay asymptotics and scalability for peer-to-peer live streaming. {\it IEEE Trans. Parallel Distrib.}, {\it 20 (10)},  1499–1511.
\item  Ferro, C.A.T. and  Segers, J. (2003). Inference for Clusters of Extreme Values. {\it J. R. Statist. Soc. B.}, {\it 65}, 545-556.
\item  Leadbetter, M.R.,  Lingren, G. and  Rootz$\acute{e}$n, H. (1983). {\it Extremes and Related Properties of Random Sequence and Processes}.
   Springer, New York.
\item Markovich, N.M. (2013). Quality Assessment of the  Packet Transport  of Peer-to-Peer Video Traffic in High-Speed Networks. {\it Performance Evaluation}, {\it 70},  28–44.
\item Markovich, N.M. (2014). Modeling clusters of extreme values. {\it Extremes}, {\it 17(1)}, 97-125.
\item Roberts, G. O.,  Rosenthal, J. S.,  Segers, J.,  Sousa B. (2006). Extremal Indices, Geometric Ergodicity of  Markov Chains and MCMC
	{\it Scandinavian Journal of Statistics}, {\bf 9}, 213-229.

\end{description}
\end{document}